\documentclass[12pt]{amsart}
\usepackage{amsmath}
\usepackage{graphicx}
\usepackage{hyperref}
\usepackage[latin1]{inputenc}
\usepackage{mathtools}
\usepackage{amsfonts}
\usepackage{amssymb}
\usepackage[T1]{fontenc}
\usepackage{amsthm}
\usepackage{fullpage}
\usepackage{enumitem}
\usepackage{longtable}

\newtheorem{theorem}{Theorem}[section]

\newtheorem{lemma}[theorem]{Lemma}
\newtheorem{proposition}[theorem]{Proposition}

\theoremstyle{definition}
\newtheorem{definition}{Definition}[section]

\theoremstyle{remark}

\numberwithin{equation}{section}

\newcommand{\Fm}{\mathbb{F}_{q^m}}

\DeclareMathOperator{\Tr}{Tr}

\title{Primitive pairs of rational functions with prescribed traces over finite fields}


\author{Shikhamoni Nath}
\address{Department of Mathematical Sciences, Tezpur University, Tezpur, Assam, 784028, India}
\email{shikha@tezu.ernet.in}

\author{Dhiren Kumar Basnet}
\address{Department of Mathematical Sciences, Tezpur University, Tezpur, Assam, 784028, India}
\email{dbasnet@tezu.ernet.in}

\begin{document}
	%\vspace{.3cm}
	\begin{abstract}
		Let $q$ be a positive integral power of some prime $p$ and $\mathbb{F}_{q^m}$ be a finite field with $q^m$ elements for some $m \in \mathbb{N}$. Here we establish a sufficient condition for the existence of a non-zero element $\epsilon \in \mathbb{F}_{q^m}$, such that $(f(\epsilon), g(\epsilon))$ is a primitive pair in $\mathbb{F}_{q^m}$ with two prescribed traces, $\Tr_{{\mathbb{F}_{q^m}}/{\mathbb{F}_q}}(\epsilon)=a$ and $\Tr_{{\mathbb{F}_{q^m}}/{\mathbb{F}_q}}(\epsilon^{-1})=b$, where $f(x), g(x) \in \mathbb{F}_{q^m}(x)$ are rational functions with some restrictions and $a, b \in \mathbb{F}_q$. Also, we show that there exists an element $\epsilon \in \mathbb{F}_{q^m}$ satisfying our desired properties in all but finitely many fields $\mathbb{F}_{q^m}$ over $\mathbb{F}_q$. We also calculate possible exceptional pairs explicitly for $m\geq 9$, when degree sums of both the rational functions are taken to be 3.
	\end{abstract}
	
	\maketitle
	
	\section{Introduction}
	
	Let $\mathbb{F}_{q^m}$  denote a field with $q^m$ elements where $q$ is some positive integral power of prime $p$ and $m\in \mathbb{N}$.	
	A  {\em primitive} element of $\mathbb{F}_{q^m}$ is a generator of the multiplicative group $\mathbb{F}^*_{q^m}$. There are $\phi (q^m-1)$ primitive elements in $\mathbb{F}_{q^m}$, where $\phi$ is the Euler's totient function. These primitive elements are of paramount importance due to their extensive use in areas like cryptographic schemes, for example in the Diffie-Hellman key exchange protocol, pseudo-random generators etc.
	\par For an element $\epsilon\in\mathbb{F}_{q^m}$, the elements $\epsilon,\epsilon^q, \epsilon^{q^2},\ldots, \epsilon^{q^{m-1}}$ are called conjugates of $\epsilon$ with respect to $\mathbb{F}_q$. Sum of all these conjugates is called the trace of $\epsilon$ over $\mathbb{F}_q$ and is denoted by $\Tr_{\mathbb{F}_{q^m}/\mathbb{F}_q}(\epsilon)$. Many researchers worked in the direction of establishing the existence of primitive elements with some prescribed conditions such as traces, which can be studied from articles like \cite{cao}, \cite{sharma1} \cite{sharma2}, among others. While majority of them focused on finding primitive pairs of the form $(\epsilon,f(\epsilon))$ for some function $f(x)$, Cao and Wang\cite{cao} attempted to find primitive functional values  $f(\epsilon)$ without considering $\epsilon$ to be primitive. But there was some error which was later addressed by Booker et al. in \cite{tru}. Here authors have finally investigated the existence of a non-zero element $\epsilon \in \mathbb{F}_{q^m}$ such that $\epsilon+\epsilon^{-1}$ is primitive and $\Tr_{\mathbb{F}_{q^m}/\mathbb{F}_q}(\epsilon)=a$, $\Tr_{\mathbb{F}_{q^m}/\mathbb{F}_q}(\epsilon^{-1})=b$ for some arbitrary member $a,b \in \mathbb{F}_q$. Moreover existence of primitive pairs where both components are rational values i.e., of the form $(f(\epsilon), g(\epsilon))$ for some rational functions $f(x), g(x) \in \mathbb{F}_{q^m}(x)$ has not yet been thoroughly investigated. So in this article we establish a sufficient condition for the existence of a non-zero element $\epsilon \in \mathbb{F}_{q^m}$, which is not necessarily primitive such that $(f(\epsilon), g(\epsilon))$ is a primitive pair in $\mathbb{F}_{q^m}$ with two prescribed traces, $\Tr_{{\mathbb{F}_{q^m}}/{\mathbb{F}_q}}(\epsilon)=a$ and $\Tr_{{\mathbb{F}_{q^m}}/{\mathbb{F}_q}}(\epsilon^{-1})=b$, where $f(x), g(x) \in \mathbb{F}_{q^m}(x)$ are rational functions with some restrictions and $a, b \in \mathbb{F}_q$.

	\par Throughout this article we are using rational functions that are not of the form $f(x) = cx^jg^d(x)$, where $j$ is any integer, $d > 1$ divides $q^m-1$ and $c \in \mathbb{F}^{*}_{q^m}$, for any rational function $g(x)\in \mathbb{F}_{q^m}(x)$. Let us denote the set of all such rational functions by $\mathcal{R}$. Let for any positive integer $m$, $w(m)$ be the
	number of distinct prime divisors of $m$ and $W(m)$ be the number of distinct square free divisors of $m$. We note that $W(m)= 2^{w(m)}$.
	\par Let us consider that the degree sum of $f(x)=m_1$ and the degree sum of $g(x)=m_2$. Moreover, let $S$ be the set containing all poles and zeros of $f(x)$ and $g(x)$ along with $0$. We also denote the set, containing the pairs $(q,m)$, for which there exists a nonzero element $\epsilon \in \mathbb{F}^*_{q^m}$, which is not necessarily primitive, such that $(f(\epsilon), g(\epsilon))$ is a primitive pair with $\Tr_{{\mathbb{F}_{q^m}}/{\mathbb{F}_q}}(\epsilon)=a$ and $\Tr_{{\mathbb{F}_{q^m}}/{\mathbb{F}_q}}(\epsilon^{-1})=b$ for some prescribed values $a,b  \in \mathbb{F}_q$ and rational functions $f,g \in \mathcal{R}$, by $\mathcal{A}_p^m$. If we can show $\mathcal{A}_p^m$ is nonempty, it will imply the existence of elements satisfying the stated conditions. If $(q,m) \notin  \mathcal{A}_p^m$, we refer to such pair as exceptional pair.
	
	\par The article is structured as follows: Section \ref{2}, provides the background information needed for the study. Section \ref{3}, presents the sufficient condition we calculated and in Section \ref{4}, we use this sufficient condition to show that 
	there exists element $\epsilon \in \mathbb{F}_{q^m}$ satisfying our desired properties in all but finitely many fields $\mathbb{F}_{q^m}$ over $\mathbb{F}_q$. Moreover considering $m_1=m_2=3$, we prove the following main result:
	
	\begin{theorem}\label{main}
		We have $(q,m)\in \mathcal{A}_p^m$ for every prime power $q$ and $m\geq 9$, except the following cases:\\
		$(i)$ $m=9$ and $q \in \{2,3,4,5,7,9,11,16\}$,\\
		$(ii)$ $m=10$ and $q \in \{ 2,3,4,5,7,8,11,23\}$, \\
		$(iii)$ $m=11$ and $q \in \{ 2,3,4\}$, \\
		$(iv)$ $m=12$ and $q \in \{ 2,3,4,5,7,8,11,16\}$, \\
		$(v)$ $q =2$ and $m \in \{13,14,15,16,18,20,24\}.$\\
		
	\end{theorem}

	\section{Preliminaries}\label{2}

	Given a finite abelian group $G$, a character $\chi$ of $G$ is defined as a homomorphism from $G$ to the multiplicative group $S^1:= \{ z \in \mathbb{C} : |z| = 1\}$. The characters of $G$ form a group under multiplication called $\mathit{dual\thinspace group}$ or $\mathit{character\thinspace group}$ of $G$ which is denoted by $\widehat{G}$ and indeed $\widehat{G} \cong G$. Let $\chi_1$ denote the trivial character of $G$ defined as $\chi_1(a)=1$, for all $a \in G$. A character of the group $G$ of order $d$ is usually denoted by $\chi_d$. We also denote all the characters of order $d$ by $(d)$.

	Corresponding to the two different abelian group structures present in a finite field $\mathbb{F}_{q^m}$ i.e., $(\mathbb{F}_{q^m}, +)$ and $(\mathbb{F}^*_{q^m},\cdot)$, two types of characters can be formed. We refer to them as additive character and multiplicative character denoted normally by $\psi$ and $\chi$ respectively. Multiplicative characters can be extended from $\mathbb{F}^*_{q^m}$ to $\mathbb{F}_{q^m}$ by the  rule
	\hspace{.1cm}  $\chi(0)=\begin{cases}
		0 \,\mbox{ if}\, \chi\neq\chi_1,\\
		1 \,\mbox{ if}\, \chi=\chi_1.
	\end{cases} $

	\begin{definition}{\textbf{($e$-free element)}}
		For any divisor $e$ of $q^m-1$, an element $\epsilon \in \mathbb{F}_{q^m}$ is called $e$-free if $\epsilon \neq \beta^d$ for any $\beta \in \mathbb{F}_{q^m}$ and for any divisor $d$ of $e$ other than $1$.   
	\end{definition}
	
	An element $\epsilon \in \mathbb{F}^*_{q^m}$ is primitive if and only if it is $(q^m-1)$-free. The characteristic function for the $e$-free elements of $\mathbb{F}^*_{q^m}$ is defined as follows:
	
	\begin{equation}\label{cfe}
		\rho_e: \mathbb{F}^*_{q^m} \to \{ 0,1\};   \epsilon \mapsto \theta(e) \sum_{d|e}\frac{\mu(d)}{\phi(d)}\sum_{(d)}\chi_{d}(\epsilon);
	\end{equation}
	where $\theta(e)=\frac{\phi(e)}{e}$ and $\phi, \mu$ denote the Euler's totient function and the Mobius function respectively.

	We shall also need characteristic function for prescribed values of trace. For each $a \in \mathbb{F}_q$ the characteristic function for elements $\epsilon \in \mathbb{F}_{q^m}$ such that $\Tr_{{\mathbb{F}_{q^m}}/{\mathbb{F}_q}}(\epsilon)=a$ is defined as follows:
	
	\begin{equation}\label{cft}
		\tau_a: \mathbb{F}_{q^m} \to \{ 0,1\}; \epsilon \mapsto \frac{1}{q}\sum_{\psi \in \widehat{\mathbb{F}_q}}\psi(\Tr_{{\mathbb{F}_{q^m}}/{\mathbb{F}_q}}(\epsilon)-a).
	\end{equation}
	
	The additive character $\tilde{\psi}(\epsilon)=e^{2{\pi}i \Tr(\epsilon)/p}$, for all $\epsilon \in \mathbb{F}_{q}$, where $\Tr$ is the absolute trace function from $\mathbb{F}_{q}$ to $\mathbb{F}_p$, is called the canonical additive character of $\mathbb{F}_{q}$ and every additive character of $\mathbb{F}_q$ is of the form $\psi_{\alpha}$ for some $\alpha \in \mathbb{F}_{q}$, where $\psi_{\alpha}(\beta)=\tilde{\psi}(\alpha\beta)$ for all $\beta \in \mathbb{F}_{q}$. If we denote $\widehat{\psi}(\epsilon)= \tilde{\psi}(Tr_{{\mathbb{F}_{q^m}}/{\mathbb{F}_q}}(\epsilon))$, this new canoninal character will be the lift of $\tilde{\psi}$. In particular, $\widehat{\psi}$ is the canonical character of $\mathbb{F}_{q^m}$. Therefore we have 
	\begin{align*}
		\tau_a(\epsilon) &=\frac{1}{q}\sum_{u \in \mathbb{F}_q}	\tilde{\psi}(\Tr_{{\mathbb{F}_{q^m}}/{\mathbb{F}_q}}(u\epsilon)-ua)\\
		&=\frac{1}{q}\sum_{u \in \mathbb{F}_q}	\widehat{\psi}(u\epsilon)\tilde{\psi}(-ua).
	\end{align*}

	The following lemmas will also be necessary for subsequent discussions.
	
	\begin{lemma} \label{lemma1}\cite{lai}
		Let $f(x)\in \mathbb{F}_{q^m}(x)$ be a rational function. Write $f(x)=\overset{k}{\underset{j=1}{\prod}}{f_j(x)}^{r_j}$, where $f_j(x) \in \mathbb{F}_{q^m}[x]$ are irreducible polynomials and $r_j $ are nonzero integers. Let $\chi$ be a nontrivial multiplicative character of $\mathbb{F}_{q^m}$ of square free order $d$ which is a divisor of $q^m-1$. Suppose that $f(x)$ is not of the form $cg(x)^d$ for any rational function $g(x)\in \mathbb{F}_{q^m}(x)$ and $c \in \mathbb{F}^*_{q^m}$. Then we have,
		
		$$\bigg|{\underset{\epsilon \in \mathbb{F}_{q^m},f(\epsilon)\neq  \infty}{\sum}\chi(f(\epsilon))}\bigg| \leq \bigg(\overset{k}{\underset{j=1}{\sum}}\deg(f_j)-1\bigg)q^{\frac{m}{2}}.$$
	\end{lemma}
	
	\begin{lemma}\label{lemma2} \cite{castro}
		Let $\chi$ and $\psi$ be two non-trivial multiplicative and additive characters of the field $\Fm$ respectively. Let $f,g$ be rational functions in $\Fm (x)$, where $f\neq \beta h^r$ and $g \neq h^p-h+\beta$, for any $h\in \Fm(x)$ and any $\beta\in \Fm$, and $r$ is the order of $\chi$. 
		Then
		$$\left | \underset{\alpha\in\mathbb{F}_{q^m} \setminus S}{\sum} \chi(f(\alpha))\psi(g(\alpha))\right|\leq [\deg(g_\infty)\, +\,  l_0 \,+\, l_1\,-\,l_2\,-\,2 ]q^{m/2},$$ 
		where $S$ denotes the set of all poles of $f$ and $g$, $g_\infty$ denotes the pole divisor of $g$, $l_0$ denotes  the number of distinct zeroes and poles of $f$ in the algebraic closure $\overline{\Fm}$ of $\Fm$, $l_1$ denotes the number of distinct poles of $g$ (including infinite pole) and $l_2$ denotes the number of finite poles of $f$, that are also zeroes or poles of $g$.
	\end{lemma}
	
	\begin{lemma}\cite{10}\label{kloos}
		If $\psi$ is a non-trivial additive character of $\mathbb{F}_{q^m}$ and $u, v \in \mathbb{F}_q$, not both zero, then the Kloosterman sum $K(\psi,u,v) =\underset{\epsilon \in \mathbb{F}^*_{q^m}}{\sum}\psi(u\epsilon+v\epsilon^{-1})$ satisfies
		$|K(\psi,u,v)|\leq 2q^{\frac{m}{2}}$.
	\end{lemma}
	
	Using these, we can now proceed to develop a sufficient condition.

	\section{A sufficient condition for elements in $\mathcal{A}_p^m$}\label{3}
	Let $e_1,e_2$ be two divisors of $q^m-1$. We denote the number of nonzero element $\epsilon \in \mathbb{F}_{q^m}$, which is not necessarily primitive, such that $f(\epsilon)$ is $e_1$-free and  $g(\epsilon)$ $e_2$-free with   
	$\Tr_{{\mathbb{F}_{q^m}}/{\mathbb{F}_q}}(\epsilon)=a$ and $\Tr_{{\mathbb{F}_{q^m}}/{\mathbb{F}_q}}(\epsilon^{-1})=b$, where $a,b\in \mathbb{F}_q$ and $f(x),g(x) \in \mathcal{R}$ by $N_{f,g,a,b}(e_1,e_2)$. Notations $m_1$, $m_2$ and $S$ remain same as defined in the introduction. Now using the characteristic functions (\ref{cfe}), and (\ref{cft}),  we get

	\begin{align*}
		N_{f,g,a,b}(e_1,e_2) =& \underset{\epsilon\in \mathbb{F}_{q^m}\textbackslash S}{\sum} \rho_{e_1}(f(\epsilon))\rho_{e_2}(g(\epsilon)) \tau_a(\epsilon)\tau_b(\epsilon^{-1})\\
		=& \frac{\theta(e_1)\theta(e_2)}{q^2} \underset{d_1|e_1,d_2|e_2}{\sum}\frac{\mu(d_1)\mu(d_2)}{\phi(d_1)\phi(d_2)} \underset{\epsilon\in \mathbb{F}_{q^m}\textbackslash S}{\sum}\\ &\underset{\chi_{d_1},\chi_{d_2}}{\sum}\chi_{d_1}(f(\epsilon))\chi_{d_2}(g(\epsilon))
		\underset{\psi,\psi' \in \widehat{\mathbb{F}_q}}{\sum}\psi(\Tr(\epsilon)-a)\psi'(\Tr(\epsilon^{-1})-b).
	\end{align*}
	Let us choose $u,v\in \mathbb{F}_q$ such that $\psi(\beta)=\tilde{\psi}(u\beta)$ and $\psi'(\beta)=\tilde{\psi}(v\beta)$ for all $\beta \in \mathbb{F}_{q^m}$. If $u,v$ run through all the elements in $\mathbb{F}_q$, then $\psi,\psi'$ run through all the characters in $ \widehat{\mathbb{F}_q}$. This transforms the above expression into the following form.
	
	\begin{align*}
		N_{f,g,a,b}(e_1,e_2)=& \frac{\theta(e_1)\theta(e_2)}{q^2} \underset{d_1|e_1,d_2|e_2}{\sum}\frac{\mu(d_1)\mu(d_2)}{\phi(d_1)\phi(d_2)} \underset{u,v \in \mathbb{F}_q}{\sum}\tilde{\psi}(-ua-vb)\underset{\chi_{d_1},\chi_{d_2}}{\sum}\chi_{f,g}(d_1,d_2,u,v);
	\end{align*}
	where,
	\begin{align*}
		\chi_{f,g}(d_1,d_2,u,v) =& \underset{\epsilon\in \mathbb{F}_{q^m}\textbackslash S}{\sum} \chi_{d_1}(f(\epsilon))\chi_{d_2}(g(\epsilon))\widehat{\psi}(u\epsilon+v\epsilon^{-1}).
	\end{align*}
	
	Our aim is to estimate a lower bound for the term $\chi_{f,g}(d_1,d_2,u,v)$, so that eventually we can evaluate a lower bound for $N_{f,g,a,b}(e_1,e_2)$. We divide the process into four sub-cases and find a lower bound of $N_{f,g,a,b}(e_1,e_2)$ in each of these cases, that leads to the formation of a general lower bound for all. We abbreviate $N_{f,g,a,b}(e_1,e_2)$ as $ N_{\text{case number}}$ for our convenience.
	\newline
	
	\textbf{Case I:} $(d_1,d_2,u,v)=(1,1,0,0)$\\
	In this case we have 
	$$\big|\chi_{f,g}(1,1,0,0)\big|\geq q^m-(m_1+m_2+1).$$
	This thereby provides,
	\begin{align*}
		N_{I} \geq \frac{\theta(e_1)\theta(e_2)}{q^2}\{ q^m-(m_1+m_2+1) \} .
	\end{align*}
	
	\textbf{Case II:} $d_1=1$  (provided $(d_2,u,v)\neq (1,0,0)$ ) \\
	Here,\\	 
	$$N_{f,g,a,b}(e_1,e_2) =  \frac{\theta(e_1)\theta(e_2)}{q^2} \underset{d_2|e_2}{\sum}\frac{\mu(d_2)}{\phi(d_2)}\underset{\substack{u,v \in \mathbb{F}_q \\ (d_2,u,v) \neq (1,0,0)}} {\sum}\tilde{\psi}(-ua-vb)\underset{\chi_{d_2}}{\sum}\chi_{f,g}(1,d_2,u,v).$$
	When $d_2=1$, we use Lemma \ref{kloos}, to arrive at $$\big|\chi_{f,g}(1,1,u,v)\big|\leq 2q^{\frac{m}{2}}.$$
	If we consider $d_2 \neq 1$, the following situations occur, which are tackled using the Lemma \ref{lemma1} and Lemma \ref{lemma2}.
	\begin{align*}
		\big|\chi_{f,g}(1,d_2,0,0)\big| & \leq (m_2-1)q^{\frac{m}{2}},\\
		\big|\chi_{f,g}(1,d_2,u,0)\big| & \leq (m_2)q^{\frac{m}{2}},\\
		\big|\chi_{f,g}(1,d_2,0,v)\big| & \leq (m_2)q^{\frac{m}{2}},\\
		\big|\chi_{f,g}(1,d_2,u,v)\big| & \leq (m_2+2)q^{\frac{m}{2}}.
	\end{align*}
	
	The number of cases when $d_1=d_2=1$ is $(q^2-1)$ and when $d_1=1, d_2 \neq 1$ is $(W(e_2)-1)q^2$. Therefore combining all these cases we get,
	
	\begin{align*}
		N_{II} &\geq -   \frac{\theta(e_1)\theta(e_2)}{q^2}q^{\frac{m}{2}}[2(q^2-1)+(W(e_2)-1)\{(q-1)^2(m_2+2)\\
		&+(q-1)(m_2)+(q-1)m_2+(m_2-1)\}]\\
		&=- \frac{\theta(e_1)\theta(e_2)}{q^2}q^{\frac{m}{2}}[W(e_2)\{q^2m_2+2q^2-4q+1\}-q^2m_2+4q-3].
	\end{align*}
	
	\textbf{Case III:} $d_2=1$  (provided $(d_1,u,v)\neq (1,0,0)$) \\
	Following the same procedure as in Case II, we can show
	\begin{align*}
		N_{III} &\geq -\frac{\theta(e_1)\theta(e_2)}{q^2}q^{\frac{m}{2}}[W(e_1)\{q^2m_1+2q^2-4q+1\}-q^2m_1+4q-3].
	\end{align*}

	\textbf{Case IV:} $d_1 \neq 1, d_2 \neq 1$\\
	In this case $\chi_{f,g}(d_1,d_2,u,v)$ can be transformed into
	\begin{align*}
		\chi_{f,g}(d_1,d_2,u,v) &= \underset{\epsilon\in \mathbb{F}_{q^m}\textbackslash S} {\sum} \chi_{q^m-1}(f^{k_1}(\epsilon)g^{k_2}(\epsilon))\widehat{\psi}(u\epsilon+v\epsilon^{-1}). 
	\end{align*}
	Here, we have considered $\chi_{d_1}=\chi^{k_1}_{q^m-1}$ and $\chi_{d_2}=\chi^{k_2}_{q^m-1}$ for some multipicative character $\chi_{q^m-1}$ of order $q^m-1$, where $k_1 \in \{ 0,1, \dots ,q^m-2 \}$ and $k_2 = \frac{q^m-1}{d_2}$.
	
	\par Then again using Lemma \ref{lemma1} and Lemma \ref{lemma2}, we derive the following,
	
	\begin{align*}
		\big|\chi_{f,g}(d_1,d_2,0,0)\big| & \leq (m_1+m_2-1)q^{\frac{m}{2}},\\
		\big|\chi_{f,g}(d_1,d_2,u,0)\big| & \leq (m_1+m_2)q^{\frac{m}{2}},\\
		\big|\chi_{f,g}(d_1,d_2,0,v)\big| & \leq (m_1+m_2)q^{\frac{m}{2}},\\
		\big|\chi_{f,g}(d_1,d_2,u,v)\big| & \leq (m_1+m_2+2)q^{\frac{m}{2}}.
	\end{align*}
	
	Now considering the number of possibilities in each case we get,
	
	\begin{align*}
		N_{IV} &\geq -\frac{\theta(e_1)\theta(e_2)}{q^2}q^{\frac{m}{2}}(W(e_1)-1)(W(e_2)-1)[(q-1)^2(m_1+m_2+2)+(q-1)\{(m_1+m_2-1)\\
		&+(m_1+m_2-1)\}+(m_1+m_2-1)]\\
		&= -\frac{\theta(e_1)\theta(e_2)}{q^2}q^{\frac{m}{2}}(W(e_1)-1)(W(e_2)-1)[q^2m_1+q^2m_2+2q^2-4q+1].
	\end{align*}

	Let us denote $T=W(e_1)W(e_2),R=W(e_1)+W(e_2)$ and sum up all these cases to find the general lower bound.
	\begin{align*}
		N_{f,g,a,b}(e_1,e_2)& =	N_{I}+	N_{II}+	N_{III}+N_{IV}\\
		N_{f,g,a,b}(e_1,e_2)&\geq \frac{\theta(e_1)\theta(e_2)}{q^2} [\{q^m-(m_1+m_2+1)\}-q^{m/2}\{W(e_1)(q^2m_1+2q^2-4q+1)\\
		&+4q-q^2m_1-3+W(e_2)(q^2m_2+2q^2-4q+1)+4q-q^2m_2-3\\
		&+(T-R+1)(q^2m_1+q^2m_2+2q^2-4q+1)\}]\\
		&\geq \frac{\theta(e_1)\theta(e_2)}{q^2} [q^m-q^{m/2}\{q^2(m_1W(e_1)+m_2W(e_2))+2q^2T+q^2(m_1+m_2)T\}].
	\end{align*}
	
	Thus we have $N_{f,g,a,b}(e_1,e_2)>0$, whenever 
	\begin{align*}
		q^{\frac{m}{2}-2}>(m_1W(e_1)+m_2W(e_2))+(m_1+m_2+2)W(e_1)W(e_2).
	\end{align*}
	
	To show existence of our desired element, we need to set $e_1=e_2=q^m-1$. This leads us to the following sufficient condition.
	\begin{align}
		q^{\frac{m}{2}-2} \geq (m_1+m_2)W(q^m-1)+(m_1+m_2+2)W(q^m-1)^2.
	\end{align} 
	Since $W(q^m-1)^2>W(q^m-1)$ and $m_1+m_2+2>m_1+m_2$, we can find another version of the above sufficient condition:
	\begin{align}\label{suff}
		q^{\frac{m}{2}-2} \geq 2(m_1+m_2+2)W(q^m-1)^2.
	\end{align}

	Sieving techniques were introduced by Cohen and Huczynska [\cite{cohen},\cite{huc}] to get improved results in such situations. Kapetanakis [\cite{kapetanakis}] also explored their findings. Following them Lemma \ref{primesievelemma1} was developed.
	\begin{lemma}[\cite{sharma2}, Lemma 3.3]\label{primesievelemma1}
		Let $d|q^m-1$ and $\{p_1,p_2,...,p_r\}$ be the collection of all primes dividing $q^m-1$ but not $d$. Then,
		\begin{align*}
			N_{f,g,a,b}(q^m-1,q^m-1)\geq \underset{i=1}{\overset{r}{\sum}}N_{f,g,a,b}(d,p_id)+\underset{i=1}{\overset{r}{\sum}}N_{f,g,a,b}(p_id,d)-(2r-1)N_{f,g,a,b}(d,d).
		\end{align*}
	\end{lemma}
	
	Aided by the proof of [\cite{sharma1}, Lemma 3.2], the next lemma follows.
	\begin{lemma}\label{primesievelemma2}
		Let $d|q^m-1$ and $p$ be a prime dividing $q^m-1$ but not $d$. Then,
		\begin{align*}
			\big|N_{f,g,a,b}(pd,d)-\theta(p)N_{f,g,a,b}(d,d)\big| &\leq \theta(p)\theta(d)^2(m_1+m_2+2)q^{\frac{m}{2}}W(d)^2,\\
			\big|N_{f,g,a,b}(d,pd)-\theta(p_i)N_{f,g,a,b}(d,d)\big| &\leq \theta(p_i)\theta(d)^2(m_1+m_2+2)q^{\frac{m}{2}}W(d)^2.
		\end{align*}
	\end{lemma}
	Using Lemma \ref{primesievelemma1} and Lemma \ref{primesievelemma2}, we deduce the following prime sieving inequality.
	\begin{theorem}\label{sieve suff}
		Let $d|q^m-1$ and $\{p_1,p_2,...,p_r\}$ be the collection of all primes dividing $q^m-1$ but not $d$. Let $l =1-2\underset{i=1}{\overset{r}{\sum}}\frac{1}{p_i}$ and $L=\frac{2r-1}{l}+2$. Assuming $l>0$, if 
		\begin{align}
			q^{m/2-2}>(m_1+m_2)W(d)+(m_1+m_2+2) W(d)^2l;
		\end{align}
		then $(q,m) \in \mathcal{A}_p^m$.
	\end{theorem}

	\section{Numerical computations}\label{4}
	In the following Proposition \ref{prop}, we will show that there exists element $\epsilon \in \mathbb{F}_{q^m}$ satisfying our desired properties in all but finitely many fields $\mathbb{F}_{q^m}$ over $\mathbb{F}_q$. We require the next lemma in order to do this.
	
	\begin{lemma}[ \cite{huc}, Lemma 3.7]\label{D}
		
		For any positive integer $M$ and positive real number $\nu$, $W(M)\leq DM^{1/\nu}$, where $D=\underset{i=1}{\overset{t}{\prod}}\frac{2}{p_i^{1/\nu}}$ and $p_1, p_2,\dots, p_t$ are the primes $\leq 2 ^\nu$ that divide $M$.
	\end{lemma}
	
	\begin{proposition}\label{prop}
		Let $q$ be a prime power and $m\geq 5$. Then, there exists an element  $\epsilon \in \mathbb{F}_{q^m}$ such that $f(\epsilon)$ and $g(\epsilon)$ both are  $(q^m-1)$-free i.e., primitive elements of $\mathbb{F}_{q^m}$ with   
		$Tr_{{\mathbb{F}_{q^m}}/{\mathbb{F}_q}}(\epsilon)=a$ and $Tr_{{\mathbb{F}_{q^m}}/{\mathbb{F}_q}}(\epsilon^{-1})=b$, where $a,b\in \mathbb{F}_q$ and $f(x),g(x) \in \mathcal{R}$ , in all but finitely many fields $\mathbb{F}_{q^m}$. 
	\end{proposition}
	
	\textbf{Proof:} From our sufficient condition (\ref{suff}) we have,
	
	$$q^{\frac{m}{2}-2} > 2(m_1+ m_2+2)W(q^m-1)^2.$$
	
	Using Lemma \ref{D} in the above inequality, we deduce the following
	\begin{align}\label{inequality1}
		q^{\frac{m}{2}-2} > 2(m_1+ m_2+2)D^2q^{2m/\nu}.
	\end{align}
	
	Now, taking logarithm on both sides, we get 
	\begin{align}\label{log}
		\log q > \frac{\log2(m_1+m_2+2)D^2}{\frac{m}{2}-2-\frac{2m}{\nu}}.
	\end{align}

	The right hand side of the above inequality is a decreasing function of $m$, if $\frac{m}{2}-\frac{2m}{\nu}>2$ i.e. $\frac{\nu-4}{2\nu}>\frac{2}{m}$, which is valid for $m\geq5$ and holds for suitable values of $\nu$.

	Therefore, there exists a number $q_\circ$ such that Inequality \ref{log} holds for $q\geq q_\circ$ and $m\geq5$.
	\par Now, for $q< q_\circ$, we consider the following equivalent Inequality \ref{inequality1},
	
	\begin{align}\label{inequality2}
		m > \frac{\log2(m_1+ m_2+2)D^2+2\log q}{(\frac{1}{2}-\frac{2}{\nu})\log q}.
	\end{align} 
	The denominator of the right hand side is positive if $\nu >4$. Hence corresponding to each $2\leq q < q_\circ$ there exists a number $m_q$ such that Inequality \ref{inequality2} holds for $m\geq m_q$.
	Therefore, based on the discussion above, we can conclude that there exists element $\epsilon \in \mathbb{F}_{q^m}$ satisfying our desired properties in all but finitely many fields $\mathbb{F}_{q^m}$ over $\mathbb{F}_q$.
	\par 	Now we proceed in a more particular direction considering $m_1=3,m_2=3$. In this case the sufficient condition takes the form 
	\begin{align}\label{new suff}
		q^{\frac{m}{2}-2} \geq 2^{2(2+w(q^m-1))}.
	\end{align} 
	\begin{lemma}[\cite{tru}, Lemma 4.1]\label{robin's bound}
		For $n\geq3$,
		\begin{align*}
			w(n)\leq \frac{1.385\log n}{\log \log n}.
		\end{align*}
	\end{lemma}
	If we set $x=q^m-1$, then using Lemma \ref{robin's bound} we can derive the following,
	\begin{align*}
		(x+1)^{\frac{1}{2}-\frac{2}{m}} \geq 2^{4+\frac{2.77 \log x}{\log \log x}}.
	\end{align*}
	Let for $m\geq m'$ , (\ref{new suff}) is true if 
	\begin{align}\label{for m'}
		(x+1)^{\frac{1}{2}-\frac{2}{m'}} \geq 2^{4+\frac{2.77 \log x}{\log \log x}}.
	\end{align}
	Now given an $m'$, we can find $x'$ for which (\ref{for m'}) holds for all $x \geq x'$. Since $q \geq 2$, we have $x \geq 2^{m'}-1$. So if the value of $x'$ is less than $2^{m'}-1$, this would imply that whenever $m\geq m'$, the pair $(q,m) \in \mathcal{A}_p^m$ for all choices of $q$. Through trial and error method we find out that $m'=109$ is a suitable choice to proceed further as for $m'=109$, we can find the value $x'\approx 3.79 \times 10^{32}$ which is less than $2^{109}-1.$ This instantly proves that  for all $m\geq 109$, we have $(q,m) \in \mathcal{A}_q^m$ for all $q$. Hence we are left with to check for $m\leq 108$.
	\par We note that, if $p_1$, $p_2$, $\dots$ , $p_i$ denote the distinct prime factors of $q^m-1$ i.e., $w(q^m-1)=i$, then 
	$$q> (p_1p_2 \cdots p_i+1)^{1/m}.$$
	Combining this observation with (\ref{new suff}) we get that $(q,m)\in \mathcal{A}_p^m$ if 
	\begin{align}\label{for i}
		(p_1p_2\cdots p_i+1)^{\frac{1}{2}-\frac{2}{m}}>2^{2(2+i)}.
	\end{align}
	For $m=108$, (\ref{for i}) holds for $i\geq22$. Hence we only need to consider $q$ such that $i \leq 21$. Since lower the value of $i$ implies higher the chances of satisfying the sufficient conditions, we consider the worst case scenario i.e., $i=21$. For $m=108$ and $i=21$, (\ref{new suff}) holds for all prime powers $q$. Repeating this same procedure we rule out the cases $m=107$, $106$, $\dots$, $97$. For $96\geq m \geq 9$, we find out the largest possible value of $i$ for which (\ref{for i}) doesn't hold. For each of those pair $m$ and $i$, we check validity of (\ref{new suff}). In each case we get a lower bound on $q$. Let us denote this lower bound as $q_{m,i}$. So for all $q \leq q_{m,i}$, we list out the pairs $(q,m)$ and calculate exact value of $i$ in each case and verify with (\ref{new suff}). The pairs that fail to satisfy that condition are further treated with sieving techniques. For each such pair we try to calculate suitable value of $ d $ so that (\ref{sieve suff}) holds. The pairs for which we could not find such $d$ are considered as possible exceptional pairs and thus we complete the proof of Theorem \ref{main}.
	
	\par List of possible exceptional pairs for $m \geq 9$ :
	$$\{ (2,9),(3,9),(4,9),(5,9),(7,9),(9,9),(11,9),(16,9),(2,10),(3,10),(4,10),(5,10),(7,10),$$
	$$
	(8,10),(11,10),(23,10),(2,11),(3,11),(4,11),(2,12),(3,12),(4,12),(4,12),(5,12),(7,12),$$
	$$(8,12),(11,12),(16,12),(2,13),(2,14),(2,15),(2,16),(2,18),(2,20),(2,24)\}$$ 
	Although our sufficient condition is valid for all $m \geq 5$, we could not check for the pairs with $5\leq m\leq8$, due to limitation of time and resources to handle such large calculations.

	\section{Declarations}
	\textbf{Conflict of interest} The authors declare no competing interests.
	
	\textbf{Ethical Approval} Not applicable.
	
	\textbf{Funding} The first author is supported by NFOBC fellowship (NBCFDC Ref. No. 231610154828).
	
	\textbf{Data availability} Not applicable.

\end{document}